\documentclass[11pt]{amsart}
\usepackage{amsmath, amsthm, amsfonts, ifpdf}
\usepackage[dvipsnames,usenames]{color}

\theoremstyle{plain}
\newtheorem{theorem}{Theorem}[section]
\newtheorem*{theorem*}{Theorem}

\newtheorem{pro}[theorem]{Proposition}
\newtheorem{Def}[theorem]{Definition}
\newtheorem{lem}[theorem]{Lemma}

\theoremstyle{definition}
\newtheorem*{Def*}{Definition}

\newtheorem{Rem}[theorem]{Remark}

\numberwithin{equation}{section}

\newcommand{\bpo}{\begin{pro}}
\newcommand{\epo}{\end{pro}}
\newcommand{\be}{\begin{equation}}
\newcommand{\ene}{\end{equation}}
\newcommand{\br}{\begin{Rem}}
\newcommand{\er}{\end{Rem}}
\newcommand{\bl}{\begin{lem}}
\newcommand{\el}{\end{lem}}
\newcommand{\bd}{\begin{Def}}
\newcommand{\ed}{\end{Def}}
\newcommand{\ben}{\begin{enumerate}}
\newcommand{\een}{\end{enumerate}}
\newcommand{\bp}{\begin{proof}}
\newcommand{\ep}{\end{proof}}
\newcommand{\beq}{\begin{equation*}}
\newcommand{\eeq}{\end{equation*}}
\newcommand{\bear}{\begin{eqnarray*}}
\newcommand{\eear}{\end{eqnarray*}}
\newcommand{\bt}{\begin{theorem}}
\newcommand{\et}{\end{theorem}}
\newcommand{\bst}{\begin{split}}
\newcommand{\est}{\end{split}}

\newcommand{\bal}{\begin{aligned}}
\newcommand{\eal}{\end{aligned}}

\renewcommand{\P}{\partial}
\newcommand{\F}[2]{\frac{#1}{#2}}
\newcommand{\la}{\langle}
\newcommand{\ra}{\rangle}

\newcommand{\R}{\mathbb{R}}

\newcommand{\bnb}{\bar{\nabla}}
\newcommand{\nb}{\nabla}

\newcommand{\RM}{Riemannian manifold}

\newcommand{\Sc}{\varepsilon}

\renewcommand{\H}{\mathbb{H}}
\newcommand{\Ta}{\Theta}

\newcommand{\vp}{\varphi}

\newcommand{\tw}{\tilde{\omega}}
\newcommand{\PLH}{{\mkern-1mu\times\mkern-1mu}}

\newcommand{\uh}{\tilde{h}}
\newcommand{\wf}{\phi}

\def\XXint#1#2#3{{\setbox0=\hbox{$#1{#2#3}{\int}$}
    \vcenter{\hbox{$#2#3$}}\kern-.5\wd0}}

\makeatletter
\def\@citestyle{\m@th\upshape\mdseries}
\def\citeform#1{{\bfseries#1}}
\def\@cite#1#2{{%
  \@citestyle[\citeform{#1}\if@tempswa, #2\fi]}}
\@ifundefined{cite }{%
  \expandafter\let\csname cite \endcsname\cite
  \edef\cite{\@nx\protect\@xp\@nx\csname cite \endcsname}%
}{}
\makeatother
\begin{document}
\title[Existence, Regularity and Asymptotic behaviors of IMCF]{Inverse mean curvature flows in warped product manifolds}
\author{Hengyu Zhou}
\address{Department of Mathematics, Sun Yat-sen University, No. 135, Xingang Xi Road, Guangzhou, 510275, People's Republic of China}
\email{zhouhy28@mail.sysu.edu.cn}
\date{\today}
\begin{abstract} We study inverse mean curvature flows of starshaped, mean convex hypersurfaces in warped product manifolds with a positive warping factor $\wf(r)$. If $\wf'(r)>0$ and $\wf''(r)\geq 0$, we show that these flows exist for all times, remain starshaped and mean convex. Plus the positivity of $\wf''(r)$ and a curvature condition we obtain a lower positive bound of mean curvature along these flows independent of the initial mean curvature. We also give a sufficient condition to extend the asymptotic behavior of these flows in Euclidean spaces into some more general warped product manifolds.
\end{abstract}
\maketitle
\section{Introduction}
   In this article we study inverse mean curvature flows in warped product manifolds. Our objective is to propose that for mean convex, starshaped hypersurfaces in warped product manifolds such flows have good evolution behaviors under certain weak conditions. \\
  \indent Fix $n\geq 3$. Throughout this paper we assume that $N$ is a $(n-1)$-dimensional closed Riemannian manifold with a smooth metric $\sigma$ and $\wf(r)$ is a smooth positive warping function on $\R^+$. A warped product manifold is the set $\{(x,r)\}$ for $x\in N, r\in \R^+$ endowed with the warped metric $dr^2+\wf^2(r)\sigma$, written as $N\PLH_{\wf}\R^+$. This general definition includes Euclidean spaces, Hyperbolic spaces \cite{GL14}, Fuchsian manifolds \cite{HLZ16}, anti-de Sitter-Schwarzchild manifolds \cite{BHW14} as well as Schwarzschild-Kottler spaces \cite{LW12} etc (also see \cite{BA08}). A hypersurface is starshaped if it is a graph over $N$ in $N\PLH_{\wf}\R^+$. \\
  \indent Assume $\Sigma$ is a mean convex closed smooth hypersurface in a Riemannian manifold $M$. An inverse mean curvature flow of $\Sigma$ is a smooth map $F:\Sigma\times [0, T]\rightarrow M$ satisfying
\be\label{def:imcf}
\F{\P F}{\P t}=\F{\vec{v}}{H}
\ene
where $F(.,0)$ is the identity on $\Sigma$. Here $H$ and $\vec{v}$ are the mean curvature and normal vector of $\Sigma_t=F(.,t)(\Sigma)$
 respectively.\\
\indent The investigation of inverse mean curvature flows are closely related to various warped product manifolds.  Gerhardt \cite{Ger90} and Urbas \cite{Urb90} first studied inverse curvature flows for admissible starshaped hypersurfaces in Euclidean spaces in which inverse mean curvature flow is a special case. Also see a recent work of Mullins \cite{Mul16}. These authors obtained long time existence and asymptotic behaviors that the flow converges uniformly to a slice of warped product manifolds after a metric scaling if $h(r)=r$ and the Ricci curvature of $N$ is positive. A similar long time existence result was established in \cite{Ger11}. However Hung-Wang \cite{HW15} showed that no such asymptotic behaviors in Hyperbolic spaces by examples. Rotationally symmetric spaces $S^n\PLH_h\R^+$ were investigated by Ding \cite{DQ11} and Scheuer \cite{SJ17}.\\
\indent On the other hand inverse mean curvature flow is a powerful tool to explore geometric properties of ambient manifolds. Some interesting applications were given by Huisken-Ilmanen \cite{HI01} for Riemannian Pernose inequalities, by Bray-Neves \cite{BN04} for the Poincar\'{e} conjecture for 3-manifolds with Yamabe invariant and by Brendle-Hung-Wang \cite{BHW14}, Kwong-Miao \cite{KM14} and Li-Wei \cite{LW12} for various geometric inequalities in certain classes of warped product manifolds. \\
\indent With delicate integral techniques Huisken-Ilmanen \cite{HI08} obtained high regularity properties of inverse mean curvature flows in Euclideans space: a sharp lower bound for the mean curvature of starshaped hypersurfaces along the flows, independent of initial mean curvature. A central tool is the Sobolev inequality in \cite{MS73}. Without this inequality  Li-Wei \cite{LW12} discovered similar estimates in Schwarzschild spaces. A consequence of high regularity properties is that an inverse mean curvature flow of bounded weakly mean convex, starshaped $C^1$ hypersurface will become smooth instantly in the settings of \cite{HI08}, \cite{LW12}. Other recent progresses on inverse mean curvature flows can be found in \cite{All15}, \cite{GLW16},\cite{CL16} and \cite{CL15} etc.\\
\indent The results mentioned above imply some nontrivial connections between inverse mean curvature flows and warped product manifolds. We are motivated to summarize such kinds of connections.  \\
\indent This paper is organized as follows. In section 2 we develop some preliminary facts for latter applications. Section 3 is devoted to show Theorem \ref{thm:LT}. Provided $ h'(r)>0$ and $h''(r)\geq 0$ we obtain the long time existence and the preservation of starshaped, mean convex properties of inverse mean curvature flows for starshaped, mean convex initial hypersurfaces in warped product manifold $N\PLH_h\R^+$. Comparing to most previous results (\cite{SJ17},\cite{LW12},\cite{Ger90},\cite{Urb90}) we discover the essential role of the condition $h'(r)>0$ and $h''(r)\geq 0$ in the evolution of inverse mean curvature flows in warped product manifolds. \\
\indent In section 4, we demonstrate a lower positive bound for mean curvature of starshaped inverse mean curvature flows in warped product manifolds under the assumption: $h'(r)>0$, $h''>0$ and $hh''-h'^2+\rho\geq 0$ for all $r>0$ where $Ric_N\geq (n-2)\rho\sigma$( see \eqref{condition:c1} and Theorem \ref{thm:lowerbound}). This bound is independent of the curvature of the initial surface and only depends on the support function $\omega=h(r)\la \vec{v},\P_r\ra$ of the initial surface, the time $t$ and the positivity of $h''(r)$. Comparing to \cite{HI08}, our proof does not rely on the Sobolev inequality \cite{MS73} which is very restrictive in general Riemannian manifolds (see \cite{HS74}). Inspired by Li-Wei \cite{LW12}, our new ingredients here are the assumption mentioned above is sufficient to obtain Theorem \ref{thm:lowerbound}. Similarly as in \cite{HI08} and \cite{LW12}, inverse mean curvatuer flows can smooth a $C^1$ starshaped hypersurface with bounded, weak nonnegative mean curvature in the setting of Theorem \ref{thm:lowerbound} (see Theorem \ref{thm:app}). \\
    \indent In section 5 we discuss the question: to what extent there eixst similar asymptotic behaviors of starshaped, mean convex inverse mean curvature flows in warped product manifolds as these in Euclidean spaces. \cite{HW15} and \cite{LW12} gave some negative and positive answers for this question respectively. In views of these existing examples, we answer this question positively in Theorem \ref{thm:finite_bound} provided that (a) $h'$ is positive with an upper uniformly bound;(b)$h''(r)\geq 0$ and satisfies $$h''(r)=O(h^{-(1+\alpha)}(r))$$ as $r\rightarrow \infty$; (c) the Ricci curvature of $N$ is positive. Here $\alpha$ is any given positive number.
    Condition (a) is necessary due to the negative answer for Hyperbolic spaces in \cite{HW15} (see Remark \ref{Rmk:sum2}). Our condition is almost optimal (see Remark \ref{Rmk:sum3}).\\
      \indent In this paper we do not consider the case that $h'(r)$ goes to infinity as $r\rightarrow \infty$. First in \cite{SJ17} Scheuer has an excellent description on inverse mean curvature flows in rationally symmetric spaces for this case. Second we are working on the evolution of inverse curvature flows in warped product manifolds in a coming paper. We hope to deal this issue with some new ingredients.\\
      \indent  Most of this work were done when the author visited Nanjing University from September 2015 to March 2016. He is very grateful for discussions with Prof. Jiaqiang Mei, Prof. Yalong Shi and Prof. Yiyan Xu. He is also very grateful for encouragements from Prof. Zheng Huang and Prof. Lixin Liu. He would like to thank the comments from Prof. Yong Wei and Prof. Mat Langford.

\section{Preliminary Facts}
In this section we develop preliminary facts needed for latter applications.
\subsection{The geometry of graphs} 
   Let $\Sigma$ be a starshaped smooth hypersurface in $N\PLH_h\R^+$, namely, it is represented by $\{(x,r(x)):x\in N\}$ where $r(x)$ is a positive function on $N$. Let $\{x_i\}$ be a local coordinate in $N$. We denote $\F{\P}{\P x_i}$ by $\P_i$ and $\F{\P}{\P r}$ by $\P_r$. The metric $\sigma$ on $N$ takes the form
\be
\sigma =\sigma_{ij}dx_i dx_j
\ene
 As in \cite{BHW14}, we define a positive function $\vp(x)$ on $N$ by
\be\label{def:vp}
\vp(x)=\Phi(r(x))
\ene
where $\Phi(r)$ is a positive function satisfying $\Phi'(r)=\F{1}{h(r)}$ . The covariant derivatives of $\vp$ are denoted by $\vp_i$ and $\vp_{ij}$. The upward normal vector of $\Sigma$ is
\be\label{def:v}
  \vec{v}=\Ta(\P_r -\F{1}{h(r)}\vp^i\P_i)\quad \text{with} \quad \Ta=\la \vec{v}, \P_r\ra = \F{1}{\sqrt{1+|D\vp|^2}}
\ene
Here $\vp^i=\sigma^{ik}\vp_k$ and $D\vp$ is the gradient of $\vp$ on $N$. $\Ta$ is called as the angle function of $\Sigma$. \\
\indent In what follows  let $X_i$ denote the vector $\P_i + r_i\P_r$. We obtain a local frame $\{X_i\}_{i=1}^n$ along $\Sigma$. With this frame the metric matrix of $\Sigma$ is
 \be\label{def:g}
   g_{ij}=\la X_i,X_j\ra =h^2(r)(\sigma_{ij}+\vp_i\vp_j)
 \ene
 with the inverse $g^{ij}=\F{1}{h^2(r)}(\sigma^{ij}-\Ta^2\vp^i\vp^j)$.\\
 \indent  The second fundamental form of $\Sigma$ is computed by
 \be\label{def:sf}
 h_{ij}=-\la \bnb_{X_i}X_{j},\vec{v}\ra
 \ene
 where $\bnb$ is the covariant derivative of $N\PLH_h\R^+$. It is well-known that $V=h(r)\P_r$ is a conformal vector field on $N\PLH_h\R^+$ (see \cite{Mon99}) such that $
 \bnb_{X}V=h'(r)X
 $ for any smooth tangent vector filed $X$.
 A straightforward computation yields that
 \begin{align}
h_{ij}&=-h(r)\Ta(\vp_{ij}-h'(r)(\sigma_{ij}+\vp_i\vp_j))\\
h_j^i &=g^{ik}h_{kj}=-\F{\Ta}{h(r)}(\tilde{\sigma}^{ik}\vp_{kj}-h'(r)\delta_{ij})\label{eq:exp:hij}
\end{align}
where $\tilde{\sigma}^{ik}=\sigma^{ik}-\Ta^2 \vp^i\vp^k$. Since $H=\sum_{i=1}^{n-1}h_i^i$, we obtain that
\bpo \label{pro:ep} The mean curvature of $\Sigma=\{(x,r(x)):x\in N\}$ in $N\PLH_h\R^+$ is
   \be\label{eq:H}
    H=-\F{\Ta}{h(r)}(\tilde{\sigma}^{ik}\vp_{ij}-(n-1)h'(r))
    \ene
    where $n-1$ is the dimension of $N$.
\epo
 Our computation is similar to those in \cite{BHW14}, \cite{Ger11} and \cite{Urb90} except replacing the sphere with any Riemannian manifold.
\subsection{Some facts on curvature tensor}  Let $\nb$ be the covariant derivative of a Riemannian manifold $N$. The Riemann curvature tensor of $N$ is given by
  \be
   R(X,Y)Z=\nb_{Y}\nb_{X}Z-\nb_{X}\nb_{Y}Z+\nb_{[X,Y]}Z
  \ene
We write $ R(X,Y,Z,W)$ for $\la R(X,Y)Z, W\ra$. Now we introduce a commuting formula for covariant derivatives and skip its proof since it is a direct verification.
\bl\label{lm:inchange} Suppose $\vp_i dx^i$ is a covector on $N$. Then
  $$\vp_{ijk}=\vp_{ikj}+R_{kjip}\vp^p $$
  where $R_{kijp}=R(\P_k,\P_j,\P_i,\P_p)$ and $\vp^p=\sigma^{pk}\vp_k$.
\el
Let $Ric_N, Ric$ denote the Ricci curvature of $N$ and $N\PLH_h\R^+$ respectively.  Their relation is given as follows.
\bpo[Proposition 9.106 in \cite{BA08}]\label{pro:Ric} Let $n-1$ be the dimension of $N$. Then
\be
 Ric=Ric_N-(h(r)h''(r)+(n-2)h'^2(r))\sigma-(n-1)\F{h''(r)}{h(r)}dr^2
\ene
\epo
\subsection{Nonparametric form} Let $\Sigma$ be a smooth hypersurface in $N\PLH_h\R^+$. The parametric form of the inverse mean curvature flow is the function $F_1:\Sigma\times [0, T)\rightarrow N\PLH_h\R^+$ as the solution in \eqref{def:imcf}.
 The nonparametric form of the inverse mean curvature flow is
 the function $F_2:\Sigma\times [0, T)\rightarrow N\PLH_h\R^+$ satisfying
 \be\label{def:non:imcf}
 (\F{\P }{\P t}F_2)^{\bot}=\F{\vec{v}}{H}
 \ene
 where $\bot$ is the projection into the normal bundle of $F_2(.,t)$. In fact $F_1(x,t)=F_2(g(x,t),t)$ where $g(x,t)$ is a family of diffeomorphisms on $\Sigma$. For its derivation, we refer to Ecker-Huisken \cite{EH89}. We will use two forms of inverse mean curvature flows interchangeably. \\
 \indent Suppose $\Sigma_t$ is a nonparametric starshaped mean convex inverse mean curvature flow with a representation $(x, r(x,t))$ in $N\PLH_h\R^+$. Then $r(x,t)$ satisfies
  \be\label{eq:r:nonparametric}
  \F{\P r}{\P t}=\F{1}{H\Ta}
\ene
In terms of $\vp(x,t)=\Phi(r(x,t))$, this is equivalent to
\be\label{eq:vp:nonparametric}
\F{\P \vp}{\P t}=\F{1}{Hh(r)\Ta}=\F{1}{H\omega}
\ene
where $\omega$ is the support function $\la h(r)\P_r,\vec{v}\ra$.
  \section{Long time existence}
  In this section we establish a long time existence of inverse mean curvature flows as follows.
   \bt\label{thm:LT}
   Let $h(r)$ be a positive function on $(0,\infty)$ satisfying
   \be\label{condition:c1}
   h'(r)>0, h''(r)\geq 0
   \ene for all $r \in (0, \infty)$. Then in the warped product manifold $N\PLH_h\R^+$ the inverse mean curvature flow \eqref{def:imcf} with any starshaped, mean convex initial hypersurface remains starshaped, mean convex and exists for all times.
 \et
  Similar long time existence results in various rotationally symmetric spaces were obtained by \cite{DQ11},\cite{LW12},\cite{BHW14} and \cite{SJ17}. Also see \cite{Ger11} and \cite{Urb90} for general curvature flows in Euclidean spaces.
 \subsection{Evolution Equations}
 Now we record some evolution equations for geometric quantities along the flows in {\RM}s. A key point is the role of Ricci curvature of ambient Riemannian manifolds in these equations.
\bpo\label{pro:mleq} Let $\Sigma_t$ be an inverse mean curvature flow in a Riemannian manifold $M$. Let $g_{ij}$ and $h_{ij}$  denote the metric and the second fundamental form on $\Sigma_t$ respectively. Then
\begin{enumerate}
\item $\P_t g_{ij}=\F{2}{H}h_{ij}$ and $\P_t g^{ij}=-\F{2}{H}g^{ik}h_{kl}g^{lj}$;
\item  $\P_t \vec{v}=\F{\nb H}{H^2}$ where $\nb$ is the gradient on $\Sigma_t$;
\item $\P_td\mu=d\mu$;
\item $\P_t H =\F{\Delta H}{H^2}-2\F{|\nb H|^2}{H^3}-\F{1}{H}(|A|^2+\bar{R}ic(\vec{v},\vec{v}))$;
\item $\P_t h_{ij}=\F{\Delta h_{ij}}{H^2}-\F{2}{H^3}\nb_i H\nb_j H-2\F{\bar{R}_{0i0j}}{H}
+(\F{|A|^2}{H^2}+\bar{R}^k_{\ 0k0})h_{ij}-\F{1}{H^2}(2h^p_l \bar{R}^l_{\ ijp}+h^p_j\bar{R}^k_{\ ikp}+h^p_i\bar{R}^k_{\ jkp})
-\F{1}{H^2}(\bnb_k\bar{R}_{0ji}^{\ \ \ k}+\bnb_{i}\bar{R}^k_{\ 0kj})$;
\end{enumerate}
Here $\bar{R}, \bar{R}ic$ are the Riemann curvature and Ricci curvature of $M$. $d\mu$ is the volume form of $\Sigma_t$.
\epo
\bp The derivation of the first four identities is very classical (see \cite{Hui86}). So we skip the detail here. As for the last one, applying the definition of $h_{ij}$ in \eqref{def:sf} and $\P_t=H^{-1}\vec{v}$ one derives that
\begin{align}
  \P_t h_{ij}&=\la \bnb_{\P_t}\bnb_{X_i}\vec{v}, X_j\ra +\la \bnb_{X_i}\vec{v},\bnb_{\P_t}X_j\ra;\notag\\
     &=\bar{R}(X_i, \P_t, \vec{v}, X_j)+\la\bnb_{X_i}\bnb_{\P_t}\vec{v}, X_j\ra+\la \bnb_{X_i}\vec{v},\bnb_{X_j}\F{\vec{v}}{H}\ra;\notag\\
     &=-\F{1}{H}\bar{R}_{0i0j}+\F{1}{H^2}\nb_{i}\nb_{j}H -\F{2}{H^3}\nb_{i}H\nb_{j}H+\F{1}{H}h_{ik}h_{jl}g^{kl}\label{id:I}
     \end{align}
In the last line we have used conclusion (2). Thus conclusion (5) follows from the expression of $\nb_i\nb_j H$, referred as the Simons' identity (see Lemma 2.1 in \cite{Hui86}).
\ep
  Following \cite{HI08} for a starshaped hypersurface $\Sigma$ in $N\PLH_{h}\R^+$ we define the support function $\omega$ and the modified speed function $u$ as follows
\be
   \omega =h(r)\la \vec{v},\P_r\ra =h(r)\Ta\quad u=\F{1}{H\omega}
\ene
where $\vec{v}$ is the upward normal vector of $\Sigma$.
\bpo\label{pro:ss} Let $\Sigma_t$ be the mean convex, starshaped inverse mean curvature flow in \eqref{def:imcf}. Let $\Delta, \nb$ be the Laplacian operator and the covariant derivative on $\Sigma_t$ respectively. Then $\omega$ and $u$ satisfy that
\begin{align}
\P_t \omega&=\F{\Delta \omega}{H^2}+\F{|A|^2}{H^2}\omega-\F{h}{H^2}Ric(\vec{v},\P_r^T);\label{eq:omega}\\
&=\F{\Delta \label{eq:2omega}\omega}{H^2}+\F{|A|^2}{H^2}\omega+\omega\F{(1-\Ta^2)}{H^2h^2}K_h\\
\P_t u&=\F{\Delta u}{H^2}-2u^{-1}\F{|\nb u|^2}{H^2}-2\F{\nb_i H\nb_i u}{H^3}-(n-1)\F{h''}{h}u^3\omega^2\label{eq:u}
\end{align}
where $\P_r^T$ is the component of $\P_r$ along $\Sigma_t$ and $K_h$ is
\be
(n-2)(hh''-h'^2)+Ric_N(v_N,v_N)
\ene
Here $\Sigma_t$ has the representation $(x,r(x,t))$, $v_N=0$ if $\Ta=1$ and $v_N=\F{Dr(x,t)}{|Dr(x,t)|}$ if $\Ta\neq 1$ where $Dr$ is the gradient of $r(x,t)$ on $N$.
\epo
\br Notice that $\sigma(v_N, v_N)=1$ if $\Ta \neq 1$. Assuming $Ric_{N}\geq (n-2)\rho \sigma$ \eqref{eq:2omega} becomes
\be\label{eq:st:omega}
\P_t\omega\geq\F{\Delta \omega}{H^2}+\F{|A|^2}{H^2}\omega+\omega\F{(n-2)}{H^2h^2}(1-\Ta^2)(hh''-h'^2+\rho)
\ene
Suppose $hh''-h'^2+\rho\geq 0$, then
\be
\P_t\omega\geq\F{\Delta \omega}{H^2}+\F{|A|^2}{H^2}\omega
\ene
  A similar assumption $hh''-h'^2+\rho>0$ appeared as the condition (H4) of  Theorem 1.1 by Brendle \cite{Bre13} regarding the rigidity of constant mean curvature hypersurfaces in warped product manifolds.
\er
\bp Since $\bnb_{X}(h(r)\P_r)=h'(r)X$ and $\P_t\vec{v}=\F{\nb H }{H^2}$ from conclusion (2) in Proposition \ref{pro:mleq}, we compute
\be\label{eq:og1}
\P_t\omega =\F{1}{H^2}\la\nb  H, h(r)\P_r\ra +\F{h'(r)}{H}
\ene
Let $\{e_i\}$ be an orthonormal frame on $\Sigma_t$. Then $h_{ik}=\la \bnb_{e_i}\vec{v}, e_k\ra$ and
\begin{align*}
\Delta\omega&=\la\bnb_{e_i}\bnb_{e_i}\vec{v}, h(r)\P_r\ra+2\la \bnb_{e_i}\vec{v}, \bnb_{e_i} h(r)\P_r\ra+\la \vec{v}, \bnb_{e_i}\bnb_{e_i}(h(r)\P_r\ra\\
&=\la \bnb_{e_i}(h_{ik}e_k),h(r)\P_r\ra +2h'(r)\la \bnb_{e_i}\vec{v}, e_i\ra +\la \vec{v}, \bnb_{e_i}(h'(r)e_i)\ra;\\
&=h_{ik,i}\la e_k, h(r)\P_r\ra-|A|^2\omega+h'(r)H;
\end{align*}
From the Codazzi identity, we have $h_{ik,i}=h_{ii,k}+\bar{R}_{0iki}$. Notice that $\P_r^T=\la \P_r, e_k\ra e_k$. Putting those facts together, we obtain
\be\label{eq:og2}
\Delta \omega=\la \nb H, h(r)\P_r\ra+h(r)Ric(\vec{v},\P_r^T)-|A|^2\omega+h'(r)H
\ene
\eqref{eq:omega} follows from the combination of \eqref{eq:og1} and \eqref{eq:og2}.\\
 \indent Next we compute $Ric(\vec{v}, \P_r^T)$. From \eqref{def:v}, we have the decomposition
\be\label{eq:dec}
\vec{v}=\Ta \P_r -\Ta\F{Dr}{h^2(r)}\quad \P_r^T = (1-\Ta^2)\P_r +\Ta^2 \F{Dr}{h^2(r)}
\ene
 With these expressions and Proposition \ref{pro:Ric} we obtain
   \begin{align}
   Ric(\vec{v}, \P_r^T)&=\Ta(1-\Ta^2)Ric(\P_r,\P_r)-\F{\Ta^3}{h^4}Ric(Dr, Dr)\notag\\
         &=-(n-1)\Ta(1-\Ta^2)\F{h''}{h}-\Ta(1-\Ta^2)\F{1}{h^2}Ric_N( v_N, v_N)\notag\\
         &+\F{\Ta(1-\Ta^2)}{h^2}(h(r)h''(r)+(n-2)h'(r)^2)\notag\\
         &=-\F{\Ta(1-\Ta^2)}{h^2}((n-2)(hh''-h'^2)+Ric_N(v_N,v_N))
         \end{align}
where $v_N=0$ if $\Ta=1$ and $v_N=\F{Dr(x,t)}{|Dr(x,t)|}$ if $\Ta\neq 1$.
  Putting these facts together, one establishes \eqref{eq:2omega}.\\
\indent  According to \eqref{eq:omega} and conclusion (4) in Proposition \ref{pro:mleq} we compute $\P_t u$ as follows.
\begin{align}
\P_t u &=-u^2(\P_t H\omega+\P_t\omega H);\notag\\
   &=-\F{u^2}{H^2}(\omega \Delta H +H\Delta\omega)+2u^2\F{|\nb H|^2}{H^3}\omega\notag\\
   &+\F{u^2}{H}(Ric(\vec{v},\vec{v})\omega
   +h(r)Ric(\vec{v},\P_r^T))\label{eq:ms}
\end{align}
Proposition \ref{pro:Ric} implies that $Ric(\P_r, X)=0$ for any $X$ satisfying $\la X,\P_r\ra=0$. By \eqref{eq:dec}, $\P_r=\Ta\vec{v}+\P_r^T$. Thus by Proposition \ref{pro:Ric} we have
 \begin{align}
 Ric(\vec{v},\vec{v})\omega+h(r)Ric(\vec{v},\P_r^T)&=h(r)Ric(\vec{v},\P_r)\notag\\
 &=h(r)\Ta Ric(\P_r,\P_r);\notag\\
 &=-(n-1)h''(r)\Ta\label{eq:ms1}
 \end{align}
 On the other hand,
 \begin{align}
 \Delta H\omega+H\Delta \omega &=\Delta(\F{1}{u})-2\nb_i H\nb_i\omega\notag\\
 &=-\F{\Delta u}{u^2}+2\F{|\nb u|^2}{u^3}+2\F{|\nb H|^2}{H}\omega+\F{2}{H u^2}\nb_i u\nb_i H\label{eq:ms2}
 \end{align}
 Plugging \eqref{eq:ms1} and \eqref{eq:ms2} into \eqref{eq:ms}, one arrives \eqref{eq:u}.
\ep
\subsection{Starshapedness} The aim of this subsection is to show that
\bt\label{thm:mid:step} Assume $h(r)$ satisfies
$h'(r)>0$ and $h''(r)\geq 0$
  for all $r>0$. Suppose $\Sigma_t$ is an inverse mean curvature flow in \eqref{def:imcf} with a starshaped, mean convex initial hypersurface on $[0, T)$ for some finite $T>0$, then $\Sigma_t$ remains mean convex and starshaped. Moreover
$$0<C_1 \leq H(p)\leq C_0$$
for $p\in \Sigma_t$. Here $C_0$ and $C_1$ are two positive constants depending on $T$ and initial data.
 \et The first step is to prove that the warping factor $h(r)$ has an exponential growth along the flow.
  \bl\label{lm:finit} Suppose $\Sigma$ is a mean convex, starshaped hypersurface with
\be\label{eq:condition:omega}
R_2 \geq \omega \geq R_1 >0;
\ene
and $h(r)$ satisfies $h'(r)>0, h''(r)\geq 0$.
   Let $\Sigma_t$ be a starshaped, mean convex inverse mean curvature flow in \eqref{def:imcf} on $[0, T)$ in $N\PLH_h\R^+$. Then \begin{enumerate}
   \item[(1)] $h(r)$ satisfies
 \be\label{pr:fin:2}
e^{\F{t}{n-1}}R_2 \geq h(r)\geq e^{\F{t}{n-1}}R_1
 \ene
 for $(x,r)\in \Sigma_t$.
 \item[(2)] for any $t\in [0, T)$ $\Sigma_t$ lies in the region
  $N\PLH h^{-1}([R_1, R_2e^{\F{T}{n-1}}])$.
 \end{enumerate}
 \el
 \br Our proof follows from Gehardt's idea \cite{Ger90}.
 \er
\bp Conclusion (2) is obvious from conclusion (1). It is sufficient to prove (1).\\
\indent Since $h(r)$ is strictly increasing, $\Ta=1$ at the points where $h(r)$ achieves its local extreme on $\Sigma$. By \eqref{eq:condition:omega} on $\Sigma$ we have $$R_1 \leq h(r)\leq R_2$$
 Suppose $\Sigma_t$ have graphical representations $(x, r(x,t))$. Define an auxiliary function $\tilde{h}$ on the cylinder $Q_T=N\times [0, T)$ as follows
 $$
 \tilde{h}(x,t)=(\ln h(r(x,t))-\F{t}{n-1}-\ln R_2)e^{-\alpha t}
 $$
 where $\alpha$ is a fixed positive constant. Then by \eqref{def:non:imcf} $r(x,t)$ satisfies
  \be
  \F{\P r}{\P t}= \F{1}{H\Ta}\quad (\Leftrightarrow \F{\P\vp}{\P t}=\F{1}{Hh\sqrt{1+|D\vp|^2}})
  \ene
 and thus
 \be\label{eq:sfirst}
   \F{\P \tilde{h}}{\P t}=-\alpha \tilde{h}+(\F{h'(r)}{Hh(r)\Ta}-\F{1}{n-1})e^{-\alpha t};
 \ene
Let $Q_{\tilde{T}}$ be the cylinder $N\times [0, \tilde{T})$ for some $\tilde{T}\in (0,T]$. Suppose the supremum of $\tilde{h}$ on $Q_{\tilde{T}}$ is obtained at a point $(x_0, t_0)$ for some $t_0\leq\tilde{T}$. Then at this point we have
 \be \label{eq:second}
 D\vp =0,\quad D^2\vp\leq 0, \quad \Ta=1, \quad \F{\P \tilde{h}}{\P t}\geq 0;
 \ene
Therefore at $(x_0, t_0)$ Proposition \ref{pro:ep} implies that
 $$
  Hh(r)\Ta\geq (n-1)h'(r);
 $$
 Together with \eqref{eq:sfirst} and \eqref{eq:second}, this yields that
 \be
 0\geq \alpha \tilde{h}-(\F{h'(r)}{Hh(r)\Ta}-\F{1}{n-1})e^{-\alpha t}\geq \alpha \tilde{h}
 \ene
 Since $\alpha$ is strictly positive and $\tilde{h}\leq 0$ at time $t=0$, then $\tilde{h}\leq 0$ on $Q_{\tilde{T}}$. Due to the arbitrariness of $\tilde{T}$, then $\tilde{h}\leq 0$ on $Q_T$. Hence $$h(r)\leq e^{\F{t}{n-1}}R_2$$ This is the second inequality in \eqref{pr:fin:2}. The first one in \eqref{pr:fin:2} follows from a similar derivation just by considering
  $$
 \tilde{h}_1(r)=(\ln h(r)-\F{t}{n-1}-\ln R_1)e^{-\alpha t}
 $$
 where $\alpha$ is a positive constant. We skip the details.
 \ep
Now we are ready to show Theorem \ref{thm:mid:step}. Although our estimate about $H$ is very rough, it is sufficient for our purpose.
\bp Let $\Sigma_t$ be the inverse mean curvature flow of $\Sigma$ in \eqref{def:imcf} existing smoothly on a time interval $[0, T)$ for $T<\infty$.\\
 \indent Suppose the conclusion in Theorem \ref{thm:mid:step} is not correct. By continuity, we can assume $\Sigma_t$ remains mean convex and starshaped on the maximal interval $[0, T')$ where $T'<T$.\\
  \indent Assume $t\in [0, T')$ for a while. By Lemma \ref{lm:finit}, $\Sigma_t$ lies in a compact domain $N\PLH h^{-1}([R_1, R_2e^{\F{T}{n-1}}])$ only depending on $T$. Thus
   \be
   |Ric(\vec{v},\vec{v})|\leq C
   \ene
   for some constant $C$ depending on $T$.    \\
   \indent Since $|A|^2\geq \F{|H|^2}{n-1}$ and $ H>0$, equation (4) in Proposition \ref{pro:mleq} becomes
\be
\P_t H^2 \leq \F{\Delta H^2}{H^2}-4\F{|\nb H|^2}{H^2}-2\F{H^2}{n-1}+C;
\ene
By the comparison principle, we have $H(p)\leq C_0$ for $C_0$ only depending on $T$ and the initial hypersurface $\Sigma$. Since $h''(r)\geq 0$, \eqref{eq:u} is simplified as
  $$
  \P_t u\leq \F{\Delta H }{H^2}-2\F{\nb_i H \nb_i u}{H^3};
  $$ for all $t\in [0, T')$. Applying the comparison principle and \eqref{pr:fin:2} we have
$$
  H R_2e^{\F{T}{n-1}} \geq Hh(r)\geq Hh(r)\Ta=u^{-1} \geq \min_{\Sigma} u^{-1};
$$
for all $t\in [0,T')$.
As a result $H\geq C_1>0$ where $C_1=R_2^{-1}e^{-\F{T}{n-1}}\min_{\Sigma} u^{-1}$. \\
\indent From conclusion (1) of Lemma \ref{lm:finit} and $H\geq C_1>0$, we observe that
   $$
   |\F{K_h}{H^2h^2}(1-\Ta^2)|\leq  C_3
   $$
where $K_h$ is from \eqref{eq:2omega} and $C_3$ is a positive constant depending on $T$ and $\Sigma$.
Thus \eqref{eq:2omega} becomes
  \be
  \P_t\omega \geq\F{\Delta \omega}{H^2} +\F{1}{n-1}\omega-C_3\omega;
  \ene
  for all $t\in [0, T')$.\\
\indent In summary $\omega \geq e^{-C_3 T}e^{\F{T}{n-1}}R_1 $ and $C_0\geq H\geq C_1>0$ for all $t\in [0, T')$. Since $T'<T$, by continuity $\Sigma_t$ is still starshaped and mean convex at time $T'$. The maximality of $T'$ gives a contradiction. Thus $T'=T$ and the starshapedness and mean convex property are preserved along $\Sigma_t$.
\ep
\subsection{The proof of Theorem \ref{thm:LT}} It is well-known that for inverse mean curvature flows, the positive lower bound of mean curvature implies its regularity. This fact was firstly shown by Smoczyk \cite{Smo00} in dimension 2, and pointed by Huisken-Ilmanen \cite{HI08} for general Riemannian manifolds. For the convenience of readers and completeness, we give its proof here.
\bl\label{lm:curin} Let $M$ be a complete Riemannian manifold. Suppose $\Sigma_t$ is a smooth mean convex solution in \eqref{def:imcf} on $[0, T)$ such that $0<H_0 \leq H \leq H_1$
 and the set $\{\Sigma_t\}_{t\in [0, T)}$ lies in a compact set $\Omega_T$ . Then the second fundamental form $A$ of $\Sigma_t$ satisfies the estimate
 \begin{align*}
    |A| &\leq C(H_0, H_1, \Omega_T)
 \end{align*}
for any $t\in [0, T)$ where $C(H_0, H_1, \Omega_T)$ is a positive constant depending on $H_0, H_1$ and $\Omega_T$.
\el
\bp As in \cite{HI08} and \cite{LW12} we define a new tensor
 $$
 \eta_i^j=Hh_i^j
 $$
  By Proposition \ref{pro:mleq} tedious computations yield the evolution equation of $\eta_i^j$ as follows.
\begin{align*}
\P_t\eta_j^i&=\F{\Delta \eta_i^j}{H^2}-2\F{\la\nb H, \nb \eta_j^i\ra}{H^3}-2\F{\nb^i H\nb_j H}{H^2}-2\bar{R}_{0i0j}-2\F{\eta_k^i\eta_j^k}{H^2}\\
&-\F{1}{H}\{\bar{R}^{li}_{\ \ lq}h^{q}_j+2\bar{R}^{li}_{\ \ jq}h^{q}_l+\bar{R}_{qlj}^{\ \ \ l}h^{qi}+\bnb^i\bar{R}_{0lj}^{\ \ \ l}+\bnb^l \bar{R}_{0j\ l}^{\ \ i}\}
\end{align*}
Since $\Omega_T$ is a compact set and $0<H_0\leq H\leq H_1<\infty $, for any $t\in [0, T)$ $\eta_j^i$ satisfies that
\be
  \P_t\eta_j^i\leq \F{\Delta \eta_i^j}{H^2}-2\F{\la\nb H, \nb \eta_j^i\ra}{H^3}-2\F{\nb^i H\nb_j H}{H^2}-2\F{\eta_k^i\eta_j^k}{H^2}+C(|\eta|+1)
\ene
where $C$ is a positive constant depending on $\Omega_T, H_0$ and $H_1$. \\
\indent Let $k_n$ be the maximum eigenvalue of $(\eta_j^i)$. Since $\Sigma_t$ are  mean convex and $k_n$ is positive, $|\eta|\leq n^2 k_n$. Notice that $(H^{-2}\nb^i H\nb_j H)$ is a nonnegative definite matrix. According to Hamilton's maximal principle about tensor fields \cite{Ham82}, $k_n$ is bounded above by $\phi$ as the positive solution of
$$
\F{d\phi}{dt}=-2\F{\phi^2}{H_1^2}+C_3(\phi+1);
$$
where $C_3$ is some different constant depending on $\Omega_T, H_0, H_1$ and $n$. If $\phi\geq \max\{1, 2C_3 H_1^{2}\}$, then $\phi$ satisfies that
 $\P_t\phi\leq-\F{\phi^2}{H_1^2}$. In this case, $\phi\leq \F{H_1^2}{t}$. Therefore
 \be
    k_n\leq \phi\leq \max\{1, 2C_3 H_1^2,\F{H_1^2}{t}\};
 \ene
 Let $\lambda_n$ be the maximal eigenvalue of $(h_i^j)$. Again because that $\Sigma_t$ are mean convex, $|A|\leq n^2\lambda_n$. Lemma \ref{lm:curin} follows since $ \eta_j^i=Hh_j^i$, $\lambda_n\leq\F{ 1}{H_0}\max\{1, 2C_3 H_1^2,\F{H_1^2}{t}\}$
\ep
 Now we conclude Theorem \ref{thm:LT} as follows.
 \bp Let $\Sigma$ be a starshaped and mean convex hypersurface in $N\PLH_h\R^+$. Suppose its inverse mean curvature flow $\Sigma_t$ exists smoothly on the maximal interval $[0, T)$ where $T$ is finite. All constants and compact sets mentioned below shall depend on $T$ and the initial surface $\Sigma$. From conclusion (2) of Lemma \ref{lm:finit}, the flow $\Sigma_t$ lies in a compact set for $t\in [0, T)$. Theorem \ref{thm:mid:step} says that $\Sigma_t$ remains mean convex and starshaped. Moreover $C_0 \geq H\geq C_1>0$ for two constants $C_1$ and $C_0$. Thus by Lemma \ref{lm:curin} all principle curvatures are bounded above by a positive finite constant. As a result \eqref{def:non:imcf} (also see \eqref{eq:r:nonparametric},\eqref{eq:vp:nonparametric}) is a uniformly parabolic equation for $t\in [0, T)$. From the regularity theory of Krylov \cite{kry87}, its solution can be extended smoothly over time $T$. This is a contradiction to the definition of $T$. Thus we conclude that the inverse mean curvature flow $\Sigma_t$ exists for all times. Moreover it preserves starshaped and mean convex properties by Theorem \ref{thm:mid:step}. The proof is complete.
 \ep

\section{High regularity when $h''(r)>0$}
In this section we study the high regularity property of inverse mean curvature flows in warped product manifolds which is first proposed in \cite{HI08}. It is a lower bound of mean curvature along starshaped inverse mean curvature flows as follows.
\bt\label{thm:lowerbound}  Suppose $h(r)$ is a smooth positive function satisfying
 \be\label{condition:c1}
  h'(r)>0\quad h''(r)>0\quad hh''-(h')^2+\rho\geq 0
 \ene
 for all $r>0$ where $Ric_N\geq (n-2)\rho \sigma$. Let $\Sigma$ be a mean convex and starshaped smooth hypersurface in $N\PLH_h\R^+$ with the property
 \be\label{r:not}
0 <  R_1\leq \omega=\la \vec{v}, h(r)\P_r\ra \leq R_2;
 \ene
  and let $\Sigma_t$ be its inverse mean curvature flow. For $t\in (0, T)$ it holds that
 \be\label{in:H}
 H\geq e^{-\F{1}{n-1}}\sqrt{h_0(n-1)}R_1 R_2^{-1}\min\{\F{1}{\sqrt{2}}t^{\F{1}{2}}, 1\}
 \ene
 Here $h_0>0$ is the infinimum of $\F{h''(r)}{h(r)}$ on the interval $h^{-1}[R_1, R_2 e^{\F{T}{n-1}}]$.
\et
\br  Similar estimates were given in \cite{HI08} and \cite{LW12}. The former depends heavily on the Michael-Simon-Sobolev inequality (see \cite{MS73}).
The latter takes advantage of the positivity of $h''(r)$ in Schwarzchild spaces. Inspired by Lemma 4.1 in \cite{LW12} we observe that the following two estimates
\begin{gather*}
 R_2 e^{\F{t}{n-1}}\geq\omega \geq e^{\F{t}{n-1}}R_1\\
  \P_t u\leq \F{\Delta u}{H^2}-2u^{-1}\F{|\nb u|^2}{H^2}-2\F{\nb_i H \nb_i u}{H^3}-(n-1)h_0u^3\omega^2
\end{gather*}
are key ingredients to establish \eqref{in:H} (see \eqref{eq:com:omega} and \eqref{eq:uin}). Both of them are guaranteed by condition \eqref{condition:c1} and \eqref{r:not}.
\er
\bp From Theorem \ref{thm:LT} the inverse mean curvature flow $\Sigma_t$ exists for all $t$, remains starshaped and mean convex. By \eqref{r:not}, $h(r)\leq R_2$ for any point $p=(x,r)\in \Sigma$. By condition \eqref{condition:c1} the conclusion (1) of Lemma \ref{lm:finit} implies that $h(r)\leq e^{\F{t}{n-1}}R_2$ on $\Sigma_t$. Thus
   $$
   \omega \leq h(r)\leq e^{\F{t}{n-1}}R_2
   $$
   on each $\Sigma_t$. On the other hand by equation \eqref{eq:st:omega} and condition \eqref{condition:c1}, we have
     $$
     \P_t\omega\geq \F{\Delta \omega}{H^2}+\F{\omega}{n-1}
     $$
     Here we use $(n-1)|A|^2\geq H^2$ and $hh''-h'^2+\rho\geq 0$ in \eqref{condition:c1}. By \eqref{r:not} the comparison principle implies $\omega \geq e^{\F{t}{n-1}}R_1$. In summary we have the estimate
      \be\label{eq:com:omega}
      R_2 e^{\F{t}{n-1}}\geq\omega \geq e^{\F{t}{n-1}}R_1
      \ene
      on any $\Sigma_t$.\\
      \indent From now on we restrict ourself to the time interval $(0, T]$.
       By the definition of $h_0$ equation \eqref{eq:u} becomes
 \be\label{eq:uin}
 \P_t u\leq \F{\Delta u}{H^2}-2u^{-1}\F{|\nb u|^2}{H^2}-2\F{\nb_i H \nb_i u}{H^3}-(n-1)h_0u^3\omega^2
 \ene
 for all $t\in (0, T]$. Similarly as in (\cite{HI08},\cite{LW12}), we define
     $$
     v=(t-t_0)^{\F{1}{2}} u;
     $$
 where $t_0\in (0, T)$ is arbitrary but fixed. Then $v(x,t_0)=0$ on $\Sigma_{t_0}$. Notice that $\F{\Delta u}{H^2}-2\F{\nb_i H \nb_i u}{H^3}=div (\F{\nb v}{H^2})$. By equation \eqref{eq:uin} the following estimate holds for $v$
 \begin{align*}
 \P_t v &\leq \F{1}{2}(t-t_0)^{-1}v +div (\F{\nb v}{H^2})-2\F{|\nb v|^2}{H^2} v^{-1}\\
 &-(n-1)h_0(t-t_0)^{-1}v^3\omega^2;
 \end{align*}
We define $v_k=\max\{v-k,0\}$ for all $k\geq 0$ and $A(k):=\{x\in \Sigma_t: v(x,t)\geq k\}$. One derives that
 \begin{align*}
 \P_t&\int_{\Sigma_t} v_k^2 d\mu_t\leq (t-t_0)^{-1}\int_{A(k)}vv_k+\int_{A(k)}v_k^2d\mu-2\int_{A_{k}}\F{|\nb v|^2}{H^2}d\mu_t\\
 &-4\int_{A(k)}v_kv^{-1}\F{|\nb v|^2}{H^2} d\mu_t
 -2(n-1)h_0(t-t_0)^{-1}\int_{A(k)}v^3v_k \omega ^2d\mu_t\\
 &\leq(t-t_0)^{-1}\int_{A(k)}vv_k+\int_{A(k)}v_k^2d\mu-2(n-1)h_0(t-t_0)^{-1}\int_{A(k)}v^3v_k\omega^2 d\mu_t
 \end{align*}
   Since $\omega\geq e^{\F{t}{n-1}}R_1$ in \eqref{eq:com:omega} and $v\geq k$ on $A(k)$, we have the inequality
 \begin{align*}
 \P_t\int_{\Sigma_t} v_k^2 d\mu_t&\leq(t-t_0)^{-1}\int_{A(k)}vv_k+\int_{A(k)}v_k^2d\mu
 \\
 &-2(n-1)h_0(t-t_0)^{-1}k^2R_1^2e^{\F{2t}{n-1}}\int_{A(k)}vv_k d\mu_t
 \end{align*}
 To make the right hand side above nonpositive,  we can choose
  \be\label{con:kk}
    k^2(t)=\F{1}{(n-1)h_0}R_1^{-2}e^{-\F{2t_0}{n-1}}\max\{t-t_0,1\}
  \ene
  Choosing $k=k(t^*)$ we have
  $$
    \P_{t}\int_{\Sigma_{t}}v_k^2d\mu_{t}\leq 0
  $$
  for $t\in [t_0, t^*]$.
  By definition $v_k(x, t_0)\equiv 0$, the above inequality implies that $v_k(x,t)\equiv 0$ for all $t\in [t_0, t^*]$. Therefore
   $
   v(x,t)\leq k(t^*)
   $ on $\Sigma_t$ for all $t\in [t_0, t^*]$.\\
  \indent We divide $t^*\in (0, T)$ into two cases. The first case is $t^*\leq 2$. We choose $t_0=\F{t^*}{2}\leq 1$. Hence $k^2(t^*)$ takes the form
  $$
  k^2(t^*)=\F{1}{(n-1)h_0}R_1^{-2}e^{-\F{2t_0}{n-1}}\leq \F{1}{(n-1)h_0}R_1^{-2}e^{-\F{2t^*}{n-1}}e^{\F{2}{n-1}}
  $$
  since $2t_0\geq 2t^*-2$ indicates that $\max\{t^*-t_0,1\}=1$.
  By the definition of $v(x,t)$, one conclude that
  $$
  \sup_{\Sigma_{t^*}}u\leq t_0^{-\F{1}{2}}k(t^*)\leq \sqrt{2}(t^*)^{-\F{1}{2}}h^{-\F{1}{2}}_0(n-1)^{-\F{1}{2}}R_1^{-1}e^{-\F{t^*}{n-1}}e^{\F{1}{n-1}}
  $$
  The second case is $t^*\geq 2$. By choosing $t_0=t^*-1$ then $k_0$ satisfies that
  $$
  k^2(t^*) =\F{1}{(n-1)h_0}R_1^{-2}e^{-2\F{t_0}{n-1}}=\F{1}{(n-1)h_0}R_1^{-2}e^{-\F{2t^*}{n-1}}e^{\F{2}{n-1}}
  $$
  Similarly as in the first case, one has that
  $$
  \sup_{\Sigma_{t^*}}u\leq k(t_*) = h^{-1}_0(n-1)^{-\F{1}{2}}R_1^{-1}e^{-\F{t^*}{n-1}}e^{\F{1}{n-1}}
  $$
  In summary for any $t\in (0, T]$ we have
   \be\label{last:step}
  \sup_{\Sigma_t}u\leq \max\{\sqrt{2}t^{-\F{1}{2}},1\}h^{-\F{1}{2}}_0(n-1)^{-\F{1}{2}}R_1^{-1}e^{-\F{t}{n-1}}e^{\F{1}{n-1}}
  \ene
The estimate in \eqref{in:H} follows from $H=\F{1}{u\omega}$ and $\omega\leq R_2 e^{\F{t}{n-1}}$ in \eqref{eq:com:omega}. The proof is complete.
 \ep
  Now we give an application of Theorem \ref{thm:lowerbound}. Our method follows closely from Theorem 2.5 in \cite{HI08} and Theorem 1.2 in \cite{LW12}. \\
 \indent Let $\Sigma$ be a $C^1$ and oriented hypersurface in a Riemannian manifold. A measurable function $H$ is called as the weak mean curvature of $\Sigma$ if it holds that
\be
\int_{\Sigma} div Xd\mu=\int_{\Sigma}H \la X,\vec{v}\ra d\mu
\ene
for any smooth vector filed $X$ with compact support. Here $\vec{v}$ is the normal vector of $\Sigma$.
 \bt\label{thm:app} Let $h(r)$ and $N$ be given by Theorem \ref{thm:lowerbound}. Let $\Sigma_0:\Sigma\rightarrow N\PLH_h\R^+ $ be a starshaped hypersurface of class $C^1$ with measurable, bounded, nonnegative weak mean curvature $H\geq 0$ and
 \be
 0< R_1\leq\omega \leq R_2;
 \ene
 for two positive constants $R_1, R_2$. Then there is a starshaped, mean convex inverse mean curvature flow $\Sigma_t$ in $N\PLH_{h}\R^+$ on $[0, \infty)$ such that  $\Sigma_t$ converges to $\Sigma_0$ uniformly in the sense of $C^0$ as $t\rightarrow 0$.
\et
First we need an approximation lemma for starshaped hypersurfaces in warped product manifolds. We follow the spirit in (Lemma 2.6, \cite{HI08}) and (Lemma 4.2, \cite{LW12}) and obtain the following much more general result in Riemannian manifold.
\bl\label{lm:app}  Let $\Sigma_0$ be a closed $C^1$ hypersurface in a Riemannian manifold $M$ with measurable, bounded weak nonnegative mean curvature. Moreover, assume $\Sigma_0$ is not minimal. \\
\indent Then $\Sigma_0$ is of class $C^{1,\beta}\cap W^{2,p}$ for all $0<\beta<1, 1\leq p<\infty$ and can be approximately locally uniformly in $C^{1,\beta}\cap W^{2, p}$ by a family of smooth hypersurfaces $\Sigma_{\Sc}$ satisfying $H>0$ where $\Sc\in (0, \Sc_0)$ for some constant $\Sc_0$.
\el
\bp Since $\Sigma_0$ is $C^1$ and the weakly nonnegative mean curvature $H$ is uniformly bounded, standard regularity results of Allard imply that $\Sigma_0$ is of class $C^{1,\beta}\bigcap  W^{2,p}$ for all $\beta, 1\leq p<\infty$. By mollification, we can pick up a sequence of smooth hypersurfaces $\Sigma_i$ in $M$ converging locally uniformly to $\Sigma_0$ in $C^{1,\beta}\bigcap W^{2,p}$. Now we consider mean curvature flows $F_i(p,\Sc)$ starting from $\Sigma_i$ in $M$ which satisfy the mean curvature flow equation
 \begin{align}
 &\F{\P F_i(p,\Sc)}{\P\Sc}=\vec{H}=-H\vec{v}\\
 & F(p,0)=p\quad p\in \Sigma_i
 \end{align}
where $\vec{v}$ is the normal vector of $F_{i,\Sc}=F_i(.,\Sc)(\Sigma_i)$.
  In a neighborhood of $\Sigma_0$, $M$ can be represented as a product manifold $\Sigma_0\times (a,b)$ if we take the Gaussian adapted coordinate. These graphs satisfy uniformly parabolic quasilinear equations with initial data in $C^{1,\beta}\cap W^{2,p}$. By the interior Schauder regularity theory, the second fundamental form satisfies
  \be\label{estimate:A}
  |A_{i,\Sc}|^2\leq \F{c}{\Sc^{\F{1}{2}-\F{\beta}{2}}}
  \ene
  on each $F_{i,\Sc}$ where $c$ is a constant uniform in $i$, depending on the regularity of $\Sigma_0$. Furthermore, by the interior estimate in \cite{EH91a} mean curvature flows $F_{i,\Sc}$ exist smoothly on a uniform time interval $[0, \Sc_0)$ independent of $i$. \\
  \indent Without loss of generality, we can assume $F_{i,\Sc}$ always stay in a compact set of $M$. Therefore the norm of sectional curvatures and its higher derivatives are uniformly bounded for all $i$ and all $\Sc\in [0,\Sc_0)$. According to the evolution equations of $H$ and $|A|^2$ for mean curvature flows in Riemannian manifolds (\cite{Hui86}), it holds that on each $F_{i,\Sc}$
             \begin{align}
             \F{\P H_{i,\Sc}}{\P\Sc}&=\Delta H_{i,\Sc}+(|A_{i,\Sc}|^2 +Ric(\vec{v},\vec{v}))H_{i,\Sc};\label{H:maximal}\\
             \F{\P |A_{i,\Sc}|^2}{\P\Sc}&\leq \Delta |A_{i,\Sc}|^2-2|\nb A_{i,\Sc}|^2 +2|A_{i,\Sc}|^4+C|A_{i,\Sc}|^2;\label{A:maximal}
             \end{align}
 where $Ric$ is the Ricci curvature of $M$ and $C$ is some positive constant depending on the compact set containing $F_i(p,\Sc)$. Considering $\int|A_{i,\Sc}|^2d\mu$, by virtue of \eqref{A:maximal} and \eqref{estimate:A} for $p\geq 1$ we have
  \begin{align*}
  \F{d}{d\Sc}\int_{F_{i,\Sc}} |A_{i,\Sc}|^p d\mu &\leq -p(p-1)\int |A_{i,\Sc}|^{p-2}|\nb A_{i,\Sc}|^2 d\mu+p\int_{F_{i,\Sc}}|A_{i,\Sc}|^2|A_{i,\Sc}|^p d\mu\\
  &+(p-2)C\int_{F_{i,\Sc}}|A_{i,\Sc}|^p d\mu\\
   &\leq (\F{c^2p}{\Sc^{1-\beta}}+C(p-2))\int_{F_{i,\Sc}} |A_{i,\Sc}|^p d\mu
  \end{align*}
 The Gronwall's lemma implies that
   \be
   \int_{F_{i,\Sc}} |A_{i,\Sc}|^p d\mu \leq \exp(\F{c^2p}{\beta}\Sc_0^\beta+C(p-2)\Sc_0)\int_{\Sigma_0}|A|^p d\mu \leq C(p)
   \ene
 Let $i$ go to infinity, we obtain a mean curvature flow $F(.,\Sc)$ which converges to $\Sigma_0$ in the sense of $C^{1,\beta}\bigcap W^{2,p}$ satisfying
  $|A_\Sc|\leq \F{C}{\Sc^{\F{1}{2}-\F{\beta}{2}}}$. Notice that the solution solution of mean curvature flow is unique in this class. Moreover corresponding mean curvature $H_\Sc\rightarrow H$ strongly in $L^p$ along $F(.,\Sc)$ as $\Sc\rightarrow 0$, $1\leq p<\infty$. Consequently, $H_{\Sc-}=\min(H_{\Sc},0)\rightarrow H_{-}=\min (H,0)$ in $L^p$. Applying \eqref{H:maximal} on the set $\{q \in F(.,\Sc): H_\Sc(q)<0\}$ we obtain
    $$\int_{F(.,\Sc)} |H_{\Sc-}|^2 d\mu \leq \exp(c\Sc_0^\beta+C\Sc_0)\int_{\Sigma_0}|H_{-}|^2 d\mu=0 $$
    since $\Sigma_0$ has nonnegative weakly mean curvature. This implies that $H_{\Sc}\geq 0$. By the strong maximum principle of parabolic equation, $\Sigma_0$ is minimal or $H_\Sc>0$ for all $\Sc>0$. The former case is excluded by our assumption. Let $\Sigma_\Sc$ denote these hypersurfaces $F(.,\Sc)(\Sigma)$. Thus the proof is complete.
\ep
We conclude Theorem \ref{thm:app} as follows.
\bp Let $\Sigma_0$ be a $C^1$ starshaped hypersurface with weakly nonnegative bounded mean curvature in warped product manifold $N\PLH_h\R^+$. Since $h'(r)>0$, Lemma 4 in \cite{LWX14} implies that any closed hypersurface is not minimal. By Lemma \ref{lm:app}, we can construct a family of mean convex hypersurfaces $\Sigma_\Sc$ for $0<\Sc<\Sc_0$ approaching $\Sigma_0$ as $\Sc\rightarrow 0$ locally uniformly in $C^{1,\beta}\bigcap W^{2,p}$. Since the convergence is in the sense of $C^1$ and $\Sigma_0$ is starshaped, we can assume $\Sigma_{\Sc}$ is starshaped, mean convex with
 $$
 0<R_1\leq \omega\leq R_2
 $$
 on each $\Sigma_\Sc$.
  Let $\Sigma_{\Sc, t}$ denote a family of inverse mean curvature flows with initial data $\Sigma_\Sc$ for $t\in [0, 1)$. Moreover, Theorem \ref{thm:lowerbound} gives the estimate
   \be\label{eq:H:lowerbound}
   H\geq \exp(-\F{1}{n-1})R_2^{-1}R_1\sqrt{(n-1)h_0}\min(\F{1}{\sqrt{2}}t^{\F{1}{2}}, 1)
   \ene
    Note that smooth hypersurfaces $\{\Sigma_\Sc\}$ have $C^1$ uniform gradients. Let $\Sc$ approach to $0$. According to \eqref{eq:H:lowerbound} after choosing a subsequence we obtain a starshaped, mean convex inverse mean curvature flow $\Sigma_t$ for $t\in (0,1)$ such that $\Sigma_t$ converges to $\Sigma_0$ in $C^0$ as $t\rightarrow 0$. From the uniqueness of inverse mean curvature flows with initial data, Theorem \ref{thm:LT} implies that we can uniquely extend $\Sigma_t$ on $(0,1)$ into $(0,\infty)$. The proof is complete.
\ep
\section{Long time behaviors when $h'(r)$ is bounded}
In this section, we study long time behaviors of inverse mean curvature flows in
warped product manifolds when $h'(r)$ becomes uniformly bounded as $r\rightarrow \infty$. We shall seek assumptions such that these flows will have similar asymptotic behaviors as those in Euclidean spaces \cite{Ger90,Urb90} and Schwarzschild spaces \cite{LW12}.
\bt\label{thm:finite_bound} Let $\alpha$ and $C$ be positive constants. Let $(N,\sigma)$ be a closed Riemannian manifold with positive Ricci curvature. Suppose $h(r)$ satisfies
\be
 C\geq h'(r)>0,\quad  C\geq h^{1+\alpha}(r)h''(r)\geq 0 \label{eq:Ric:condition2}
\ene
 Then for any smooth, mean convex hypersurface $\Sigma$ in $N\PLH_h\R^+$, its inverse mean curvature flow $\Sigma_t$ exists smoothly for all times, remains mean convex and starshaped. Moreover as $t\rightarrow \infty$
\begin{enumerate}
  \item let the flow $\Sigma_t$ be written as  the graphs $(x, r(x,t))$ in $N\PLH_h\R^+$. The rescaling function $e^{-\F{t}{n-1}}h(r(x,t))$ converges uniformly to a positive constant $c$ with an exponential decay speed;
  \item the modified height function $\vp$ in \eqref{def:vp} has the property that $|D\vp|\leq Ce^{-\beta t}$ and $|D^2\vp|\leq Ce^{-\beta t}$ for some $\beta>0$.
  \item the second fundamental form $(h_i^j)$ of $\Sigma_t$ has the property that
    $$
     |h^j_i-\F{h'(r)}{h(r)}\delta_{ij}|=\F{h'(r)}{h(r)}O(e^{-\beta t})
    $$
\end{enumerate}
\et
\br\label{Rmk:sum} Our result generalizes the following asymptotic results of inverse mean curvatures for starshaped, mean convex hypersurfaces in special warped product manifolds. Let $n\geq 3$.
\begin{enumerate}
\item The inverse mean curvature flow in $n$ dimensional Euclidean spaces obtained by Gerhardt \cite{Ger90} and Urbas \cite{Urb90} where $h(r)=r$ and $N$ is the unit sphere $S^{n-1}$.
\item $N\PLH_r\R^+$ studied by Mullins \cite{Mul16} where $h(r)=r$ and $N$ is a closed manifold with positive Ricci curvature.
\item  Schwarzchild spaces investigated by Li-Wei \cite{LW12} where $h(r)$ is a positive function satisfying  $h'(r)=\sqrt{1-2m h^{2-n}}$ and $N$ is the unit sphere $S^{n-1}$ with $m>0$. Note that $h''(r)=-(2-n)mh^{1-n} $ satisfying $h''(r)=O(\F{1}{h^2(r)})$.
\end{enumerate}
\er
\br\label{Rmk:sum2}  Conclusion (1) and (2) indicate that the sequence of scaling  metrics $|\Sigma_t|^{-\F{2}{n-1}}g$ on $\Sigma_t$ converges to the standard metric $\sigma$ on $N$. Here $|\Sigma_t|$ is the area of $\Sigma_t$. In \cite{HW15} Hung-Wang construct a counterexample such that the limiting scaling metric along an inverse mean curvature flow in $\H^3$ is not a round metric on the sphere. In their case $h(r)=\sinh(r)$ which means $h'(r)$ is unbounded.
\er
\br\label{Rmk:sum3} We argue that the assumptions in Theorem \ref{thm:finite_bound} is almost optimal as follows. Due to Remark \ref{Rmk:sum2}, $0<h'(r)\leq C$ is optimal. Moreover it is impossible find a positive function $h(r)$ with the property $ hh''(r)\geq c>0$ for some constant $c$. Otherwise because $h(r)=O(r)$, then $h''(r)\geq \F{c}{r}$. This implies that $h'(r)=O(\ln(r))$ and gives a contradiction to $h'(r)\leq C$. Therefore, we can think $h''(r)=O(\F{1}{h^{(1+\alpha)}(r)})$ is almost optimal since $\alpha$ is arbitrarily positive. To obtain the exponential decay of $|D\vp|$, we require that $Ric_N(v_N, v_N)+hh''$ at least has a positive lower bound in \eqref{eq:twoomega} for large $r$ and $v_N\neq 0$. This is only possible in the case that the Ricci curvature of $N$ is positive. Because as mentioned above $hh''(r)$ could be arbitrarily small on the interval $[r_0,\infty)$ for any $r_0>0$.
\er
\subsection{An evolution equation}
Next we investigate the evolution of $|D\vp|^2$ along inverse mean curvature flows. For abbreviation we define
 \begin{align}
 F&=Hh\Ta=\Ta^2((n-1)h'(r)-\tilde{\sigma}^{ij}\vp_{ij})\label{de:F}\\
 G_k &=(F\vp^k+\Ta^2\sigma^{ik}\vp^j\vp_{ij}-\Ta^4\vp^k\vp^i\vp^j\vp_{ij})
 \end{align}
 where $\vp^i=\sigma^{ik}\vp_{k}, \vp^k_i=\sigma^{kl}\vp_{li}$ and $\tilde{\sigma}^{ij}=(\sigma^{ij}-\Ta^2\vp^i\vp^j)$. A direct computation yields that $\F{\P F}{\P\vp_k}=-2\Ta^2 G_k$. In what follows $\tw$ denotes $\F{1}{2}|D\vp|^2$.
 \bpo\label{pro:def} Let $\Sigma_t$ be a starshaped, mean convex inverse mean curvature flow in $N\PLH_h\R^+$ with the representation $(x,r(x,t))$. Then $\tw$ satisfies
 \begin{align}\label{eq:twoomega}
 \F{\P \tw}{\P t}&= \F{1}{H^2h^2}\{\tilde{\sigma}^{ij}\tw_{ij}+ 2G_i\tw_i-\tilde{\sigma}^{ij}\vp_{ki}\vp^{k}_j\\
 &-2\tw(Ric_N(v_N, v_N)+(n-1)h(r)h''(r))\}\notag
 \end{align}
 where $v_N=0$ if $\tw =0$ or $v_N=\F{Dr(x,t)}{|Dr(x,t)|}$ if $\tw\neq 0$. Here $Ric_N$ is the Ricci curvature of $N$.
 \epo
 \br Notice that the term $\tilde{\sigma}^{ij}\vp_{ki}\vp^{k}_j$ is always nonnegative.
 \er
 \bp For $r(x,t)$ satisfies that $\F{\P r(x,t)}{\P t}=\F{1}{H\Ta}$ (see \eqref{def:non:imcf}), $\F{\P\vp}{\P t}=\F{1}{F}$. We differentiate $\tw$ with respect to $t$. Then
           \begin{align}
           \F{\P\tw}{\P t}&=\vp^k (\F{\P \vp}{\P t})_k
           = -\F{1}{F^2}\vp^k (F)_k  \\
      &=\F{\Ta^2}{F^2}(\tilde{\sigma}^{ij}\vp_{ijk}\vp^k + 2G_i\vp_{ik}\vp^k-(n-1)h(r)h''(r)\vp_k\vp^k)\label{eq:mis}
 \end{align}
 We note that
   \begin{align}
  \tilde{\sigma}^{ij}\omega_{ij}&=\tilde{\sigma}^{ij}\vp_{kij}\vp^k+\tilde{\sigma}^{ij}\vp_{ki}\vp^k_{j}\notag\\
               &=\tilde{\sigma}^{ij}\vp_{ikj}\vp^k+\tilde{\sigma}^{ij}\vp_{ki}\vp^k_{j}\notag\\
               &=\tilde{\sigma}^{ij}\vp_{ijk}\vp^k+\tilde{\sigma}^{ij}R_{jkip}\vp^p\vp^k+\tilde{\sigma}^{ij}\vp_{ki}\vp^k_{j}\label{eq:desdg}
   \end{align}
   where the covariant derivatives are taken with respect to $(N,\sigma)$. In the last line we apply Lemma \ref{lm:inchange}. Moreover,
   \begin{align}
 \tilde{\sigma}^{ij}R_{jkip}\vp^p\vp^k&=(\sigma^{ij}-\Ta^2\vp^i\vp^j)R_{jkiq}\vp^q\vp^k;\notag\\
                      &=\sigma^{ij}R_{jkiq}\vp^q\vp^k-\Ta^2\vp^i\vp^j R_{jkiq}\vp^q\vp^k;\notag\\
                      &=Ric_N(\P_k,\P_q)\vp^q\vp^k\notag \\
                      &= Ric_N(D\vp, D\vp)\label{eq:desdf}
   \end{align}
   where the third line is from $\vp^i\vp^j R_{ikjq}\vp^q\vp^k=0$. Here $R$ is the Riemann curvature tensor of $N$. Note that
   $Ric_N(D\vp, D\vp)= 2\tw Ric(v_N, v_N)$ if $\Ta\neq 1$. The proposition follows from combining \eqref{eq:mis},\eqref{eq:desdf} with \eqref{eq:desdg}.
 \ep
   \subsection{The proof of Theorem \ref{thm:finite_bound}} By Theorem \ref{thm:LT} we can assume that $\Sigma_t$ is a starshaped, mean convex inverse mean curvature flow in $N\PLH_h\R^+$ on $[0,\infty)$. All functions are evaluated along this flow. Our proof is divided into three steps.
    \subsubsection{The first step} Since the modified speed function $u=\F{1}{H\omega}$, then \eqref{eq:u} also can be written as
   \be\label{eq:u:newversion}
    \P_t u=\F{\Delta u}{H^2}-2u^{-1}\F{|\nb u|^2}{H^2}-2\F{\nb_i H\nb_i u}{H^3}-(n-1)\F{h''h}{h^2H^2}u\\
    \ene
    The following result is derived from condition \eqref{eq:Ric:condition2}.
    \bl Under the setting of Theorem \ref{thm:finite_bound}, the mean curvature of $\Sigma_t$ satisfies
   \be\label{cl:bound:H}
   C_2\geq H(p)h(r)\geq C_1>0
    \ene
    for two positive constants $C_1$ and $C_2$ independent of $t$. Here $p$ denotes any point $(x,r)$ in $\Sigma_t$.
    \el
    \bp In the sequel, all computations are carried on $\Sigma_t$. By \eqref{eq:Ric:condition2}, $h''(r)\geq 0$. Applying the maximum principle in \eqref{eq:u:newversion}, we have $\max_{\Sigma_t} u\leq \max_{\Sigma} u$. This implies that
        \be\label{eq:fact:a}
           Hh(r)\geq H h(r)\Ta\geq C_1:=\F{1}{\max_{\Sigma} u}>0
        \ene
    on any $\Sigma_t$.  On the other hand we have $h(r)=O(e^{\F{t}{n-1}})$ for all $t$ by conclusion (1) in Lemma \ref{lm:finit}. Since $h''h^{1+\alpha}(r)\leq C$, we get
      \be\label{eq:fact:b}
        h''h \leq C_{0}e^{-\F{\alpha t}{n-1}}
      \ene
      for a positive constant $C_{0}$ independent of $t$. Putting \eqref{eq:fact:a} and \eqref{eq:fact:b} into \eqref{eq:u:newversion}, we obtain that
    \be\label{eq:u:newversion}
    \P_t u\geq \F{\Delta u}{H^2}-2u^{-1}\F{|\nb u|^2}{H^2}-2\F{\nb_i H\nb_i u}{H^3}-\F{C_0}{C_1^2}e^{-\F{\alpha t}{n-1}}u\\
    \ene
    Note that $\alpha$ is a positive fixed number. The maximum principle tells us that $\min_{\Sigma_t}u$ is bounded below by the solution of differential ordinary equation
     \be
     \P_t q(t)=-\F{C_0}{C_1^2}e^{-\F{\alpha t}{n-1}}q
     \ene
     with $q(0)=\min_{\Sigma}u>0$. It is easy to see that
       $$
       u(t)=\F{1}{Hh(r)\Ta}\geq q(t)\geq C_1^*>0
       $$
       for a positive constant $C_1^*$. \\
    \indent  Then $Hh(r)\leq C_2$ provided $\Ta$ has a uniformly lower positive bound. Since $\Ta=\F{1}{\sqrt{1+|D\vp|^2}}$, it is equivalent to show that $|D\vp|^2$ is uniformly bounded. Because $Ric_N\geq 0$ and $h''(r)\geq 0$, this is achieved by applying the maximal principle into \eqref{eq:twoomega}. We complete the proof.
    \ep
     \subsubsection{The second step} The second step is to show that
    \be\label{cl:vp}
    |D\vp|^2=O(e^{-2\beta t})
     \ene for some positive constant $\beta$. Since $h''\geq 0$ and $\tilde{\sigma}_{ij}\vp_{ik}\vp_{j}^k\geq 0$, by \eqref{eq:twoomega} $\tilde{\omega}=\F{1}{2}|D\vp|^2$ satisfies
  \be\label{tw:midstep}
  \F{\P \tw}{\P t}\leq  \F{1}{H^2h^2}(\tilde{\sigma}^{ij}\tw_{ij}+ 2G_i\tw_i)-2\F{\tw}{H^2h^2} Ric_N(v_N, v_N)
  \ene
  Let $s(t)$ denote the maximum of $\tilde{\omega}$ on $\Sigma_s$ for $s\in [0, t)$. Without loss of generality we assume that $s(t)>0$. Then there exist $t_0\leq t$ and $p_0\in \Sigma_{t_0}$ such that
    $$
    \tw(p_0)=s(t)
    $$
    Moreover, at this point $v_N\neq 0$ since $\Ta\neq 1$ by $s(t)>0$. By condition \eqref{eq:Ric:condition2}, $Ric_N(v_N, v_N)\geq \rho_0$ for some uniformly positive constant $\rho_0$.
    According to \eqref{tw:midstep} and $Hh\leq C_2$, the comparison principle implies that
    $$
    \P_t s(t)\leq  -2\F{\rho_0}{C_2^2}s(t)
    $$
    Let $\beta$ be the positive number $\F{\rho_0}{C_2^2}$. This implies our claim in \eqref{cl:vp}.
\subsubsection{The last step}  We consider an auxiliary function
 $$
 \tilde{h}(r)=h(r)e^{-\F{1}{n-1}t}
 $$
 Then $\uh$ satisfies
 \be\label{eq:exp:F}
 \P_t\uh=\F{h'}{H\Ta}e^{-\F{t}{n-1}}-\F{1}{n-1}\uh
 \ene
 where we use the nonparametric form of inverse mean curvature flows $\F{\P r}{\P t}=\F{1}{H\Ta}$. Moreover let $\tilde{F}$ be the right hand side of \eqref{eq:exp:F}. By the conclusion (1) of Lemma \ref{lm:finit}, $\uh$ is uniformly bounded. According to \eqref{cl:bound:H} we obtain 
   $$
   \F{h'}{H\Ta}e^{-\F{t}{n-1}}=\F{h'}{Hh\Ta}he^{-\F{t}{n-1}}
   $$
   is uniformly bounded. Here we also use $|h'|\leq C$ and $\F{1}{\Ta}=\sqrt{1+|D\vp|^2}\leq C$ by \eqref{cl:vp}. Therefore $\P_t \uh$ is uniformly bounded. On the other hand
     \be\label{mid:step:ij}
      \uh_{ij}=(h''h^2+hh'^2)e^{-\F{t}{n-1}}\vp_i\vp_j+hh'e^{-\F{1}{n-1}}\vp_{ij}
     \ene
     From the expression of $H$ in \eqref{eq:H} we derive that
     $$
     \F{\P\tilde{F}}{\P \uh_{ij}}=\F{1}{H^2h^2}(\sigma^{ij}-\F{\vp^i\vp^j}{1+|D\vp|^2})\F{1}{hh'e^{-\F{1}{n-1}}}
     $$
   By \eqref{cl:bound:H} and \eqref{cl:vp}, $(\F{\P\tilde{F}}{\P \uh_{ij}})$ is strictly definite and bounded from above and below. Therefore \eqref{eq:exp:F} is a strictly parabolic equation. Considering the gradient of $\tilde{h}$, we see that
      $$
      D\uh =h'(r)D\vp h(r)e^{-\F{t}{n-1}}=O(e^{-\beta t})
      $$
      decays with an exponential speed.\\
      \indent From the regularity theory of Krylov (\S 5.5 in \cite{kry87}), higher derivatives of $\uh$ are uniformly bounded. Applying the interpolation theorem (Lemma 6.1, \cite{Ger11} or page 371, \cite{Urb90}) upon $D^2\uh$, we obtain that $D^2\uh=O(e^{-\alpha t})$. So are all higher derivatives of $\tilde{h}$. Thus we obtain conclusion (1).\\
      \indent In view of \eqref{mid:step:ij}, we have $D^2\vp=O(e^{-\beta t})$ because $h''h\leq C$ from condition \eqref{eq:Ric:condition2} and $he^{-\F{t}{n-1}}\geq C_1>0$ in \eqref{cl:bound:H}.\\
      %
      \indent Due to the expression of $h_i^j$ in \eqref{eq:exp:hij}, we have
      \be
      |h^{i}_j-\F{h'}{h}\delta_i^j|\leq\F{h'}{h}(O(D^2\vp)+(\Ta-1))
      \leq\F{h'}{h}O(e^{-\beta t})
      \ene
  We obtain conclusion (2). The proof is complete.
  \section{Some Remarks} 
  \begin{enumerate}
  	\item In this version, we update the reference of Brian Allen's work \cite{All15}. 
  	\item  The author was informed that \cite{Tsui17} Prof. Maopei Tsui also got the same answer as Theorem 5.1 in the case of rotationally symmetric spaces when the author visited the National Taiwan University after this paper was published online.
  	\end{enumerate}
\bibliographystyle{abbrv}	
\bibliography{Ref_Thesis}
\end{document}